\documentstyle[12pt]{amsart}
\numberwithin{equation}{section}
\newtheorem{dfn}{Definition}[section]
\newtheorem{prop}[dfn]{Proposition}
\newtheorem{thm}[dfn]{Theorem}
\newtheorem{lem}[dfn]{Lemma}

\begin{document}
\title{On a $q$-analogue of the multiple gamma functions}
\author{Michitomo Nishizawa}
\address{Department of Mathematics, School of Science and Engineering,
Waseda university, 3-4-1 Ohkubo,Shinjuku,Tokyo,Japan}
\email{64m503@@cfi.waseda.ac.jp}
\maketitle


\begin{abstract}
 A $q$-analogue of the multiple gamma functions is introduced,
and is shown to satisfy the generalized Bohr-Morellup theorem.
Furthermore we give some expressions of these function.
\end{abstract}\newpage

\section{Introduction}
 In 1904, E.W.Barnes defined the multiple gamma functions
relevant to his zeta functions \cite{barnes1}, \cite{barnes2},
\cite{barnes3}, \cite{barnes4}. His works were influenced
by the study of the integral functions which are motivated
by number theory.  M.F.Vigne\'ras investigated the multiple
gamma functions on her paper \cite{vign}. She showed
these functions to satisfy the "generalized Bohr-Morellup
theorem". These functions have been studied by many authors
relevant to the determinant of Laplacian on Riemann surfaces
and the Selberg zeta function in number theory and physics
(For example \cite{vardi}, \cite{voros}).
 On the other hand, the $q$-gamma function was defined
by F.H.Jackson \cite{jackson1}, \cite{jackson2}. R.Askey showed
these function satisfies a $q$-analogue of the Bohr-Morellup
theorem \cite{askey}.
 Recently P.G.O.Freund and A.V.Zabrodin constructed a hierarchy of
S-matrices of integrable models by using the multiple gamma
function and its $q$-analogue \cite{freud}. But their $q$-analogue
remains some ambiguities because there is not the definition of the
normalization factor of these functions in their paper.
 Motivated by Vigne\'ras' and Askey's works, we construct
a $q$-analogue of the multiple gamma functions exactly, which is
a generalization of the $q$-gamma function and satisfies a $q$-analogue
of the generalized Bohr-Morellup theorem. We can see these functions
have some expressions like the gamma function. Furthermore, they are related
to a $q$-multiple zeta functions like the case of Barnes' multiple
gamma function and his multiple zeta functions \cite{ueno}, \cite{un}.

The author expresses his deep gratitude Professor Kimio Ueno for
helpful discussions and constant encouragement.

\section{A survey of the multiple gamma functions }

In this section we give a brief survey of the basic facts concerning
the gamma function and its generalization.

The gamma function are characterized Bohr-Morellup theorem.

\begin{thm}[Bohr-Mollerup]  The gamma function satisfies
  \begin{enumerate}
    \item $\Gamma(z+1)=z \Gamma(z)$,
    \item $\Gamma(1)=1$,
    \item $\frac{d^2}{dz^2}\log\Gamma(z+1)\geq0
        \mbox{ for } z\geq0,$
 \end{enumerate}
and the function satisfying (1),(2),(3) is determined uniquely.
\end{thm}

We use expression of this function as follows.

\begin{equation}
    \Gamma(z+1)=\lim_{N\to\infty}
        \frac{N!}{(z+1)(z+2)\cdots(z+N)} (N+1)^{z}.
    \label{eqn:gauss}
\end{equation}
\begin{equation}
    \Gamma(z+1)=\prod_{n=1}^{\infty}
       \left\{
        \left(1+\frac{1}{n}\right)^z
        \left(1+\frac{z}{n}\right)^{-1}
        \right\}.
    \label{eqn:euler}
\end{equation}
\begin{equation}
\Gamma(z+1)=e^{-\gamma x}
   \prod_{n=1}^{\infty}
     \left\{
     \left(1+\frac{z}{n}\right)^{-1}
     e^{\frac{z}{n}}
     \right\},
\end{equation}
where  $\gamma$ is Euler constant.\\

As generalization of the gamma function, M.F.Vign\'{e}ras
constructed the multiple gamma functions by using
a generalization of the Bohr-Morellup theorem \cite{vign}.
First, we remark a theorem due to Dufresnoy and Pisot \cite{dufr}.

\begin{thm}[Dufresnoy and Pisot]
Let $g(z)$ be a k times differentiable function and
$g^{(k)}(z)\to0 \mbox{ as }z\to\infty$, then the function $f(z)$
which satisfies $$f(z+1)-f(z)=g(z)$$ exists.  It is unique if $f(0)$ is
given. Furthermore, $f(z)$ is k times differentiable and $f^{(k)}(z)$
is increasing for $x\geq 0$
\label{thm:dufr}\end{thm}

By this theorem, if we put
    $$f_{0}(z)= log(z+1),$$
then we can determine $f_{r}(z)$ to satisfy
    $$f_{r}(z)-f_{r}(z-1)=f_{r-1}(z),\qquad f_{r}(0)=0,$$
because $f_{r}(z)$ is $(r+1)$-times differentiable and
$f^{(r+1)}_{r}(z) \to 0 \mbox{ as } z \to \infty.$
Hence we can define
    $$G_{r} (z+1)=\exp f_{r}(z),$$
then the next theorem follows (cf.\cite{vign}).

\begin{thm}
The functions $G_{r} (z)$ satisfies
\begin{equation}
\begin{aligned}
    & (1) \quad G_{r}(z+1)=G_{r-1}(z) G_{r}(z), \\
    & (2) \quad G_{r}(1)=1, \\
    & (3) \quad \frac{d^{r+1}}{dz^{r+1}}\log G_{r}(z+1)\geq0
       \quad \mbox{ for } \quad z\geq0,\\
    & (4) \quad G_{0}(z) =z.
\end{aligned} \label{eqn:mg}
\end{equation}
And meromorphic function satisfies above properties is
determined uniquely.
\label{thm:mg}\end{thm}

These functions are called ''the multiple gamma functions''. For example,
$G_{1}(z)$ is $\Gamma(z)$.

\section{The $q$-gamma function}
In this section we review the $q$-analogue of the gamma function
(The $q$-gamma function), which is known as
\begin{equation}
    \Gamma (z+1;q)=(1-q)^{-z}
        \prod_{n=1}^{\infty}\left(
            \frac{1-q^{z+n}}{1-q^{n}}
            \right)^{-1}.\label{eqn:qgamm}
\end{equation}

R.Askey showed this function satisfies a $q$-analogue of the
Bohr-Morellup theorem.

\begin{thm}[Askey]
For $0<q<1$, the function $\Gamma(z;q)$ satisfies
    \begin{enumerate}
        \item $\Gamma (z+1;q)=[z] \Gamma(z;q),$
        \item $\Gamma (1;q)=1,$
        \item $\frac{d^{2}}{dz^{2}}\log\Gamma(z+1;q)\geq0
            \mbox{ for } z\geq0.$
    \end{enumerate}
\label{thm:qBM}\end{thm}

The uniqueness of this function follows from this theorem and
theorem \ref{thm:dufr}. This theorem suggests that Vigne\`ras
method can be applied to a $q$-analogue of the multiple gamma
functions.

The $q$-gamma function has the expressions which correspond to
(\ref{eqn:gauss}), (\ref{eqn:euler}). It is easy to show

\begin{equation}
    \Gamma (z+1;q)=\lim_{N\to\infty}
        \frac{[1][2]\cdots [N]}
            {[z+1][z+2]\cdots [z+N]}
            [N+1]^{z}.
\end{equation}
and
\begin{equation}
    \Gamma (z+1;q)=\prod_{n=1}^{\infty}
        \left\{
        \left(
            \frac{[n+1]}{[n]}
        \right)^{z}
        \left(
            \frac{[z+n]}{[n]}
        \right)^{-1}
        \right\},
\end{equation}

where we use a notation
    $$[z]=\frac{1-q^{z}}{1-q}.$$

\section{A $q$-analogue of the multiple gamma functions}
In this section, we define a $q$-analogue of the multiple gamma functions
which satisfy a $q$-analogue of the generalized Bohr-Morellup theorem,
and we derive the expressions corresponding to (\ref{eqn:gauss}),
(\ref{eqn:euler}). We assume $0<q<1$.

\begin{dfn} Let $z$ be in the right half plane $\{s \in \bold C |\Re s > 0\}$
and $r \in  {\bold Z}_{\geq 0}$, we define
    $$G_{0}(z+1;q):= [z+1],$$
    $$G_{r}(z+1;q):= (1-q)^{- \binom{z}{r}} \prod_{n=1}^{\infty}
        \left\{
        \left(
        \frac{1-q^{z+n}}{1-q^{n}}
        \right)^{(-1)^{r} \binom{n+r-2}{r-1}}
        (1-q^{n})^{g_{r}(z,n)}
        \right\},\quad \mbox{for}\quad r \geq 1$$
where
    $$g_{1}(z,n):=0,\qquad g_{r} (z,n) := \sum_{m=1}^{r-1}
        (-1)^{m-1} \binom{z}{r-m} \binom{n+m-2}{m-1}\quad\mbox{for}
      \quad r\geq2$$
\end{dfn}

For example, $G_{1}(z;q)$ is $\Gamma(z;q)$.
The infinite products of these functions are absolutely
convergent, First, we prove a $q$-analogue of Theorem \ref{thm:mg}

\begin{thm}
If $\Re z > 0$, then $G_{r}(z:q)$ satisfy
    \begin{enumerate}
        \item $G_{r}(z+1;q)=G_{r-1}(z;q)G_{r}(z;q),$
        \item $G_{r}(1;q)=1,$
        \item $\frac{d^{r+1}}{dz^{r+1}}\log G_{r+1}(z+1;q)\geq0
            \mbox{ for } z\geq0,$
        \item $G_{0}(z;q)=[z].$
    \end{enumerate}
The functions satisfying such properties are determined uniquely.
\label{thm:gqBM}\end{thm}

\begin{pf}
First we remark next formulas.
\begin{lem}
\begin{enumerate}
    \item $g_{r}(0,n)=0.$
    \item $g_{r}(z,n)-g_{r}(z-1.n)=g_{r-1}(z-1,n)
        -(-1)^{r-1}\binom{n+r-3}{r-2}.$
\end{enumerate}
\end{lem}

These formula can be proved by direct calculation. Next we prove
the claim of theorem.
When $r=1$, the claims are Theorem \ref{thm:qBM}. So it is sufficient
to show the case $r \geq 2$.
(2),(4) can be proved easily. So we prove (1) and (3).

(1)We can see
{\allowdisplaybreaks
\begin{align*}
    G_{r} & (z+1;q)\\
        & = (1-q)^{- \binom{z-1}{r}-\binom{z-1}{r-1}}
            \prod_{n=1}^{\infty} \left\{
                (1-q^{z+n-1})^{(-1)^{r}\binom{n+r-2}{r-1}
                -(-1)^{r} \binom{n+r-3}{r-1}}
            \right. \\
        & \qquad \left.
            \times (1-q^{n})^{-(-1)^{r}\binom{n+r-2}{r-1}
            +g_{r}(z-1,n)+g_{r-1}(z-1,n)-(-1)^{r-1}
            \binom{n+r-3}{r-2}}
        \right\} \\
        & = (1-q)^{- \binom{z-1}{r-1}}
            \prod_{n=1}^{\infty}
            \left\{
            \left(
            \frac{1-q^{z+n-1}}{1-q^{n}}
            \right)
           ^{(-1)^{r-1}\binom{n+r-3}{r-2}}
            (1-q^{n})^{g_{r-1}(z-1,n)}
            \right\} \\
        & \qquad \times (1-q)^{- \binom{z-1}{r}}
            \prod_{n=1}^{\infty}
            \left\{
            \left(
            \frac{1-q^{z+n-1}}{1-q^{n}}
            \right)
            ^{(-1)^{r} \binom{n+r-2}{r-1}}
            (1-q^{n})^{g_{r}(z-1,n)}
            \right\} \\
        & = G_{r-1}(z;q)G_{r}(z;q).
\end{align*}}

(3) We can see
{\allowdisplaybreaks
\begin{align*}
    \frac{d^{r+1}}{dz^{r+1}} & \log G_{r}(z+1;q) \\
        & = (-1)^{r+1} \frac{d^{r+1}}{dz^{r+1}}
            \left\{
            \sum_{n=1}^{\infty}  \sum_{k=1}^{\infty}
            \binom{n+r-2}{r-1}
            \frac{q^{(z+n)k}}{k}
            \right\} \\
        & = (- \log q)^{r+1} \sum_{n=1}^{\infty}
            \sum_{k=1}^{\infty}
            \binom{n+r-2}{r-1}k^{r}q^{(z+n)k} \\
        & \geq 0,
\end{align*}}
because $\log q < 0$.\\
 Finally, we prove the uniqueness of these functions. Because of the
formula in the proof (3),
$$\frac{d^{r+1}}{dz^{r+1}}  \log G_{r}(z+1;q) \to 0$$
as $z\to \infty$. Hence the uniqueness follows from
Theorem \ref{thm:dufr}.
\end{pf}

By Theorem \ref{thm:gqBM} (1),  $G_{r}(z;q)$ can be defined
when $z \ne l + m \log q / 2 \pi $ ($l$ ; negative integer,
$m \in{\bold Z}$). Next we consider expression like the formulas
(\ref{eqn:gauss}), (\ref{eqn:euler}).

\begin{prop}
If $\Re z > 0$, then
\begin{align*}
    & (1) \quad G_{r}(z+1;q) = \lim_{N \to \infty}
        \left\{
        \frac{G_{r-1}(1;q) \cdots G_{r-1}(N;q)}
            {G_{r-1}(z+1;q) \cdots G_{r-1}(z+N;q)}
            \prod_{m=1}^{r} G_{r-m}(N+1;q)
            ^{\binom{z}{m}}
        \right\}. \\
    & (2) \quad G_{r}(z+1;q) =
        \prod_{n=1}^{\infty}
        \left\{
        \frac{G_{r-1}(n;q)}{G_{r-1}(z+n;q)}
        \prod_{m=1}^{r}
        \left(
        \frac{G_{r-m}(n+1;q)}{G_{r-m}(n;q)}
        \right)^{\binom{z}{m}}
        \right\}.
\end{align*}\label{prop:prod}
\end{prop}
\begin{pf}

(1) First we prove next formulas
\begin{lem}
(1) For $r \geq 1$
$$G_{r}(N+1;q)=(1-q)^{- \binom{N}{r}}
        \prod_{n=1}^{N}(1-q^{n})^{\binom{N-n}{r-1}} .$$
(2) For $r \geq 2$, $N \geq 1$,
    $$\sum_{m=1}^{r-1} \binom{z}{m} \binom{N}{r-m}
        = \sum_{n=1}^{N} \left\{
            \binom{z+n-1}{r-1} - \binom{n-1}{r-1}
            \right\}.$$
(3) For $N \geq 1$, $n \geq 1$, $r \geq 2$,
    $$\sum_{m=1}^{r+1} \binom{N-n}{r-m} \binom{z}{m}
        - \sum_{k=1}^{N}
        \left\{
            g_{r}(z+k-1,n)- g_{r}(k-1,n)
        \right\}
        = g_{r+1}(z,n).$$
\label{lem:qmg}\end{lem}

They are shown by induction.
We prove (1) of Proposition \ref{prop:prod}.
By Lemma \ref{lem:qmg} (1),
{\allowdisplaybreaks
\begin{align*}
&\lim_{ N \to \infty }
    \left\{
    \frac{ G_{r-1}(1;q) \cdots  G_{r-1}(N;q) }
        { G_{r-1}(z+1;q) \cdots G_{r-1}(z+N;q) }
            \prod_{m=1}^{r} G_{r-m} (N+1;q)^{\binom{z}{m}}
    \right\}  \\
    & = \lim_{N \to \infty}
        \left\{
        (1-q)^{ \sum_{k=1}^{N}
        \left\{
        \binom{z+k-1}{r-1} - \binom{k-1}{r-1}
        \right\}
        - \sum_{m=1}^{r-1} \binom{N}{r-m} \binom{z}{m}
        - \binom{z}{r}}
        \right.
    \times \prod_{k=1}^{N} \prod_{n=1}^{\infty}
        \left(
        \frac{1-q^{z+n+k-1}}{1-q^{n+k-1}}
        \right)
        ^{(-1)^{r} \binom{n+r-3}{r-2}} \\
    & \qquad \times \prod_{n=1}^{N}
        (1-q^{n})^{\sum_{m=1}^{r} \binom{N-n}{r-m}\binom{z}{m}
            - \sum_{k=1}^{N}
            \left\{
            g_{r-1}(z+k-1,n) - g_{r-1}(k-1,n)
            \right\}}\\
    & \qquad
        \left.\times
        \prod_{n=N+1}^{\infty}
        (1-q^{n})^{ - \sum_{k=1}^{N}
            \left\{
            g_{r-1}(z+k-1,n) - g_{r-1}(k-1,n)
            \right\}}
        \right\}.
\end{align*}}
By Lemma \ref{lem:qmg} (2), (3)
and
    $$\prod_{k=0}^{\infty} \prod_{n=1}^{\infty}
        \left(
       \frac{1-q^{z+n+k-1}}{1-q^{n+k-1}}
        \right)
        ^{(-1)^{r}\binom{k+r-3}{r-2}}
        = \prod_{n=1}^{\infty}
        \left(
        \frac{1-q^{z+n}}{1-q^{n}}
        \right)
        ^{(-1)^{r}\binom{n+r-2}{r-1}},$$
we obtain
\begin{eqnarray*}
& &\lim_{N \to \infty}
        \left\{
        \frac{G_{r-1}(1;q) \cdots G_{r-1}(N;q)}
            {G_{r-1}(z+1;q) \cdots G_{r-1}(z+N;q)}
        \prod_{m=1}^{r} G_{r-1}(N+1;q)^{\binom{z}{m}}
        \right\} \\
    & & \quad= \left\{
        (1-q)^{-\binom{z}{r}} \prod_{n=1}^{\infty}
        \left\{
        \left(
        \frac{1-q^{z+n}}{1-q^{n}}
        \right)^{(-1)^{r} \binom{n+r-2}{r-1}}
        (1-q^{n})^{g_{r}(z,n)}
       \right\}
        \right\}\\
     & &\qquad  \times \lim_{N \to \infty}
        \left\{
        \prod_{n=N+1}^{\infty}
            (1-q^{n})^{-\sum_{m=1}^{r}\binom{N-n}{r-m}\binom{z}{m}}
        \right\}.
\end{eqnarray*}
So, it is sufficient to show
    $$\prod_{n=N+1}^{\infty}(1-q^{n})^{\binom{N-n}{r-m} \binom{z}{m}}
        \quad \to 1 \qquad \mbox{as} \quad N \to \infty.$$
We can see
\begin{align*}
    & \left|
    \log
    \left\{
        \prod_{n=N+1}^{\infty}
        (1-q^{n})^{\binom{N-n}{r-m} \binom{z}{m}}
    \right\}
    \right|
     \leq \left|\binom{z}{m}\right|
        \sum_{n=N+1}^{\infty}
        \left|\ \binom{N-n}{r-m} \right|\
        \left|\ \log(1-q^{n}) \right|\ \\
    & \leq \frac{\left|\binom{z}{m}\right|}{(r-m)!(1-q)}
        \sum_{n=N+1}^{\infty}
        (n-N)(n-N+1) \cdots (n-N+r-m-1) q^{n} \\
    & = q^{N} \left\{
        \frac{\left| \binom{z}{m} \right|}{(r-m)!(1-q)}
        \sum_{n=1}^{\infty}
        (n+r-m-1)(n+r-m-2) \cdots (n+1).n
        q^{n}
        \right\},
\end{align*}
where we take the principal value of logarithms.
This tend to $0$ as $N \to \infty$ on any region, since
    $$\frac{1}{(r-m)!(1-q)} \sum_{n=1}^{\infty}
        (n+r-m-1)(n+r-m-2) \cdots (n+1).n. q^{n}$$
takes finite value. Thus,
    $$\left|
    \log
    \left\{
        \prod_{n=N+1}^{\infty}
        (1-q^{n})^{\binom{N-n}{r-m} \binom{z}{m}}
    \right\}
   \right| \to 0 \qquad \mbox{as} \quad N \to \infty.$$
Hence the claim follows. \\

(2) By using (1), we can prove the claim easily.
\end{pf}

\end{document}